
\input amstex

\documentstyle{amsppt}
\magnification = \magstep1
\define\sw#1{{\sb{(#1)}}}
\define\C{\bold C}
\define\cq{{\bold C}\sp 2\sb q}
\define\glq{GL\sb q(2, {\bold C})}
\define\sq{{S\sp 2\sb q(\mu ,\nu)}}
\NoBlackBoxes
\topmatter
\title
Quantum Homogeneous Spaces \\ as  Quantum Quotient Spaces
\endtitle
\author
Tomasz Brzezi\'nski
\endauthor
\affil
Institute of Mathematics, University of \L\'od\'z \\
ul. Banacha 22, 90-238 \L\'od\'z, Poland
\endaffil
\address
Institute of Mathematics, University of \L\'od\'z,
ul. Banacha 22, 90-238 \L\'od\'z, Poland
\endaddress
\email tbrzez\@ ulb.ac.be or t.brzezinski\@ damtp.cam.ac.uk\endemail
\curraddr University
of Cambridge, DAMTP, Cambridge CB3 9EW, UK. (after 1st October 1995)
\endcurraddr
\date August 1995 \enddate
\thanks Most of this paper was written during my stay at
Universite Libre de Bruxelles supported  by the European Union
Human Capital and Mobility grant. This work is
also supported  by the
grant  KBN 2 P302 21706 p01
\endthanks
\abstract
We show that certain embeddable homogeneous
spaces of a quantum group that do not correspond to
a quantum subgroup still have the structure of quantum
quotient spaces.
We propose a construction of quantum fibre bundles
on such spaces. The quantum
plane and the general quantum two-spheres are
discussed in detail.
\endabstract
\endtopmatter
\document

\baselineskip=18pt
\head 0. Introduction \endhead
A homogeneous space $X$ of a Lie group $G$ may be
always identified with the quotient space $G/G\sb 0$,
where $G\sb 0$ is a Lie subgroup of $G$. When the
notion of a homogeneous space is generalised to the case
of quantum groups or non-commutative Hopf algebras
the situation becomes much more complicated. A general
quantum homogeneous space of a quantum group
$H$ need not be a quotient space of $H$ by its
quantum subgroup. By a quantum subgroup of $H$ we
mean a Hopf algebra $H\sb 0$ such that there is
a Hopf algebra epimorphism $\pi : H\to H\sb 0$.
The quotient space is then understood as a
 subalgebra of $H$ of all points that are
fixed under the coaction of $H\sb 0$ on $H$
induced by $\pi$. A quantum homogeneous space $B$
of $H$ might be such a quotient space but it
is not in general. There is, however, a certain
class of quantum homogeneous spaces,
of which the quantum two sphere of Podle\'s
\cite{P1} is the most prominent example, that
not being quotient spaces by a quantum
subgroup of $H$, may be embedded in $H$. One
terms such homogeneous spaces {\it embeddable}
\cite{P2}. The general quantum two sphere
$\sq$ is such an embeddable homogeneous
space of the quantum group $SU\sb q(2)$,
and it is a quantum quotient space in the
above sense when $\nu = 0$.
In the latter case the corresponding subgroup
of $SU\sb q(2)$ may be identified with the
algebra of functions on $U(1)$. In this paper
we show that certain embeddable quantum homogeneous
spaces, and the general quantum two sphere
$\sq$ among them, can still be understood as
quotient spaces or fixed point subalgebras.
Precisely we show that there is a coalgebra
$C$ and a coalgebra epimorphism $\pi :H\to C$ such
that the fixed point subspace of $H$ under the
coaction of $C$ on $H$ induced from the coproduct in
$H$ by a pushout by $\pi$ is a subalgebra
of $H$ isomorphic to $B$.

The interpretation of embeddable quantum homogeneous
spaces as quantum quotient spaces allows
one to develop the quantum group gauge
theory of such spaces following the lines of
\cite{BM}. The study of such a gauge
theory becomes even more important once the
appearance of the quantum homogeneous spaces
in the A. Connes geometric description of
the standard model was annonunced \cite{C}.
For this purpose, however, one needs
to generalise the notion of a quantum principal bundle
of \cite{BM} so that a Hopf algebra
playing the role of a quantum structure group there
may be replaced by a coalgebra. We
propose such a generalisation. Since the theory of quantum principal
bundles is strictly related to the theory of Hopf-Galois
extensions (cf. \cite{S}), we thus propose
a generalisation of such extensions.

The paper is organised as follows. In Section~1
we describe the notation we use in the sequel.
In Section~2 we show a fixed point subalgebra
structure of embeddable quantum homogeneous
spaces. Next we propose a suitable generalisation
of the notion of a quantum principal bundle in
Section~3. Sections 4 and 5 are devoted to
careful study of two examples of quantum embeddable
spaces, namely the quantum plane $\cq$ \cite{M}
and the quantum sphere $\sq$ \cite{P1}.
\head 1. Preliminaries \endhead
In the sequel all the vector spaces are over the field
$k$ of characteristic not 2. $C$ denotes a
coalgebra with the coproduct $\Delta: C\to C\otimes C$ and the
counit $\epsilon: C\to k$ which satisfy the standard axioms,
cf. \cite{Sw}. For the coproduct we use the
Sweedler sigma notation
$$
\Delta c = c\sw 1\otimes c\sw 2,\qquad (\Delta
\otimes id)\circ\Delta c = c\sw 1\otimes c\sw 2
\otimes c\sw 3,\quad \text{etc.},
$$
where $c\in C$, and the summation sign and the
indices are suppressed. A vector space
$A$ is a left $C$-comodule if there exists
a map $\Delta\sb L : A\to C\otimes A$, such that
$(\Delta\otimes id)\circ\Delta\sb L =
(id\otimes\Delta\sb L)\circ\Delta\sb L$, and
$(\epsilon\otimes id)\circ\Delta\sb L = id$.
For $\Delta\sb L$ we use the explicit notation
$$
\Delta\sb L a = a\sw 1\otimes a\sw\infty ,
$$
where $a\in A$ and all $a\sw 1\in C$ and all $a\sw \infty
\in A$.

Similarly we say that a vector space
$A$ is a right $C$-comodule if there exists
a map $\Delta\sb R : A\to A\otimes C$, such that
$(\Delta\sb R\otimes id)\circ\Delta\sb R =
(id\otimes\Delta)\circ\Delta\sb R$, and
$(id\otimes \epsilon)\circ\Delta\sb R = id$.
For $\Delta\sb R$ we use the explicit notation
$$
\Delta\sb R a = a\sw 0\otimes a\sw 1 ,
$$
where $a\in A$ and all $a\sw 1\in C$ and all $a\sw 0
\in A$.

$H$ denotes a Hopf algebra with product $m: H\otimes H
\to H$, unit $1$, coproduct $\Delta: H\to H\otimes H$,
counit $\epsilon: H\to k$
and antipode $S: H\to H$. We use Sweedler's sigma notation
as before. Similarly as for a coalgebra we can
define right and left $H$-comodules. For a
right $H$-comodule $A$ we denote by $A\sp{coH}$
a vector subspace of $A$ of all elements $a\in A$
such that $\Delta\sb R a = a\otimes 1$. We say that
a right (resp. left) $H$-comodule $A$ is a
right (resp. left) $H$-comodule algebra if $A$
is an algebra and $\Delta\sb R$ (resp. $\Delta\sb L$)
is an algebra map.

A vector subspace $J$ of $H$ such that $\epsilon(J) =0$
and $\Delta J\subset J\otimes H\oplus H\otimes J$ is called
a {\it coideal} in $H$. If $J$ is a coideal in $H$ then
$C=H/J$ is a coalgebra with a coproduct $\Delta$
given by $\Delta = (\pi\otimes\pi)\circ\Delta$, where
$\pi :H\to C$ is a canonical surjection. The counit
$\epsilon$ in $C$ is defined by the commutative
diagram
$$
\CD
H      @>\epsilon>>    k\\
@V{\pi}VV            @VV{id}V\\
C      @>\epsilon>>    k
\endCD
$$
\head 2. Quantum homogeneous spaces\endhead
In this section we show that if an embeddable
quantum homogeneous space satisfies certain additional
assumption it may be identified with a quantum
quotient space.

\definition{Definition 2.1. [P2]} Let $H$ be
a Hopf algebra and $B$ be a left $H$-comodule
algebra with the coaction $\Delta\sb L : B\to
H\otimes B$. We say that $B$ is an {\it embeddable quantum
homogeneous space} or simply an {\it embeddable
H-space} if there exists an algebra inclusion
$i: B\hookrightarrow H$ such that $\Delta\circ i
= (id\otimes i)\circ\Delta\sb L$, i.e., $i$ is an
intertwiner.\enddefinition
\proclaim{Proposition 2.2}
\roster
\item A left $H$-comodule algebra $B$ is an embeddable
$H$-space if and only if there exists an algebra
character $\kappa :B\to k$ such that the linear map
$i\sb\kappa :B\to H$, $i\sb\kappa: b\mapsto b\sw 1
\kappa(b\sw\infty)$ is injective.
\item If $B$ is an embeddable $H$-space then the linear
map $\chi\sb L: B\otimes B\to H\otimes B$,
$\chi\sb L: b\otimes b' \mapsto b\sw 1\otimes b\sw\infty
b'$ is injective.
\endroster\endproclaim
\demo{Proof}
(1) If $B$ is an embeddable quantum homogeneous
space then $\kappa = \epsilon\circ i$ is a character
of $B$. Since $i$ is an intertwiner, for any
$b\in B$ we compute
$$
i\sb\kappa(b) = b\sw 1\epsilon(i(b\sw\infty)) = i(b)\sw 1
\epsilon(i(b)\sw 2) = i(b),
$$
thus $i\sb\kappa$ is an inclusion.

Conversly assume that there is a character $\kappa :B\to k$
such that $i\sb\kappa$ is injective. Then clearly $i\sb\kappa$
is an algebra inclusion. Furthermore
$$
\Delta(i\sb\kappa(b)) = b\sw 1\otimes b\sw 2\kappa(b\sw\infty)
= b\sw 1 \otimes i\sb\kappa(b\sw\infty) = (id \otimes
i\sb\kappa)\circ\Delta\sb L(b).
$$
Therefore $i\sb\kappa$ is an intertwiner as required.

(2) The canonical map $can: H\otimes H\to H\otimes H$,
$can: u\otimes v\mapsto u\sw 1\otimes u\sw 2v$ is a
linear isomorphism. Consider the diagram
$$
\CD
{} @. 0  @.  0 \\
@.    @VVV   @VVV \\
{} @. B\otimes B @>{\chi\sb L}>>  H\otimes B\\
@.    @V{i\otimes i}VV              @VV{id\otimes i}V \\
0 @>>> H\otimes H @>{can}>> H\otimes H
\endCD\tag{2-1}
$$
Clearly, both the rows and the columns of the diagram \thetag{2-1}
are
exact. Moreover for any $b, b'\in B$
$$\split
(id\otimes i)\circ\chi\sb L(b\otimes b')
   &=b\sw 1\otimes i(b\sw\infty b') = b\sw 1\otimes i(b\sw\infty)
       i(b') \\
   &=i(b)\sw 1\otimes i(b)\sw 2i(b') = can (i(b)\otimes i(b')),
\endsplit$$
and hence the diagram \thetag{2-1} is also commutative. Therefore
we conclude that the sequence $0 @>>> B\otimes B @>{\chi\sb L}>>
H\otimes B$ is exact, i.e. the map $\chi\sb L$ is injective.
\qed
\enddemo
\remark{Remark \rom{2.3}} The second assertion of
 Proposition~2.2., i.e., the injectiveness
of $\chi\sb L$, is a dual version of the statement
that the action of a group on its homogeneous
space is transitive.
\endremark
\proclaim{Proposition 2.4} Let $B$ be an embeddable $H$-space
corresponding to the character $\kappa : B\to k$. Define
a right ideal $J\sb\kappa\subset H$ by $J\sb\kappa=
\{ \sum\sb j (i\sb\kappa(b\sb j) -\kappa(b\sb j))u\sb j;
\forall b\sb j\in B, \;\;\forall u\sb j\in H\}$.
Then $J\sb\kappa$ is an $H$-coideal.\endproclaim
\demo{Proof} Clearly
$$\epsilon (i\sb\kappa(b) -\kappa(b)) =\epsilon(b\sw 1)
\kappa(b\sw\infty) -\kappa(b) = \kappa(b) -\kappa(b) =0.
$$
Furthermore
$$\split
\Delta(i\sb\kappa(b) -\kappa(b))
   &=i\sb\kappa(b)\sw 1\otimes i\sb\kappa(b)\sw 2 - \kappa(b)
       1\otimes 1\\
   &= b\sw 1\otimes (i\sb\kappa(b\sw\infty) - \kappa(b\sw\infty))
       + b\sw 1\kappa(b\sw\infty)\otimes 1 -\kappa(b)1\otimes 1 \\
   &= b\sw 1\otimes (i\sb\kappa(b\sw\infty)-\kappa(b\sw\infty))
       + (i\sb\kappa(b) -\kappa(b))\otimes 1 .
\endsplit$$
Therefore for any $b\in B$,
$$
\Delta(i\sb\kappa(b) -\kappa(b)) \in
       H\otimes J\sb\kappa \oplus J\sb\kappa\otimes H ,
$$
so that $J\sb\kappa$ is a coideal as stated. \qed\enddemo

Since $J\sb\kappa$ is a coideal of $H$, the vector
space $C = H/J\sb\kappa$ is a coalgebra and the
canonical projection $\pi : H\to C$ is a coalgebra
map. This in turn implies that $H$ is a right
$C$-comodule with the coaction $\Delta\sb R
= (id\otimes\pi)\circ\Delta :H\to H\otimes C$.
Let $H\sp{coC} = \{u\in H; \Delta\sb R u = u\otimes \pi(1)\}$.
\proclaim{Proposition 2.5} Let $B$ be an embeddable
$H$-space corresponding to the character
$\kappa: B\to k$, $J\sb\kappa$ be as
in Proposition~2.4 and $C=H/J\sb\kappa$. Then:
\roster
\item $H\sp{coC}$ is a subalgebra of $H$.
\item $B$ is a subalgebra of $H\sp{coC}$.
\endroster\endproclaim
\demo{Proof} (1) Since $\ker\pi = J\sb\kappa$
is a right ideal in $H$ there is a natural right
action  $\rho\sb 0: C\otimes H \to C$ of $H$ on $C$
given by the commutative diagram
$$
\CD
H\otimes H @>m>>  H \\
@V{\pi\otimes id}VV  @VV{\pi}V \\
C\otimes H @>{\rho\sb 0}>> C
\endCD
$$
In other words for any $a\in C$ and $u\in H$,
$\rho\sb 0(a,u) = \pi(vu)$, where $v\in\pi\sp{-1}(a)$.
For any $u, v\in H\sp{coC}$ we compute
$$\split
\Delta\sb R (uv)
   &= u\sw 1 v\sw 1 \otimes \pi(u\sw 2 v\sw 2)\\
   &= u\sw 1 v\sw 1 \otimes \rho\sb 0 (\pi(u\sw 2), v\sw 2) \\
   &= uv\sw 1\otimes \rho\sb 0(\pi(1), v\sw 2) \\
   &= uv\sw 1\otimes \pi(v\sw 2) = uv\otimes \pi(1).
\endsplit$$
Therefore $uv\in H\sp{coC}$ and $H\sp{coC}$ is a subalgebra
of $H$ as required.

(2) For any $b\in B$ we compute
$$\split
\Delta\sb R(i\sb\kappa(b))
   &=i\sb\kappa(b)\sw 1\otimes \pi(i\sb\kappa(b)\sw 2)
     =b\sw 1\otimes \pi(i\sb\kappa(b\sw\infty))\\
   &=b\sw 1\otimes \kappa(b\sw\infty)\pi(1)
     =i\sb\kappa(b)\otimes\pi(1).
\endsplit$$
Hence $i\sb\kappa :B\hookrightarrow H\sp{coC}$ is the
required algebra inclusion. \qed\enddemo

Proposition~2.5. shows therefore that if $H\sp{coC}
\subset i\sb\kappa(B)$ then the embeddable $H$-space
$B$ may be identified with the quantum quotient
space $H\sp{coC}$.
\head 3. A possible generalisation of quantum principal
bundles \endhead
Once an $H$-embeddable space $B$ is identified
with a quotient space $H\sp{coC}$, it is natural to
view $H$ as a total space of a principal bundle
over $B$. Therefore one would like to apply the
general theory of quantum principal bundles of [BM]
to this case too. In general, however, neither
$C$ is a Hopf algebra nor, if it happens
to be a Hopf algebra, $C$ is a quantum subgroup
of $H$. Hence the induced coaction of
$C$ on $H$ is not an algebra map. Therefore to
develop a gauge theory on embeddable homogeneous
spaces one needs to generalise the theory of
quantum principal bundles. In this section we
propose such a generalisation. It is based on
a simple observation that the structure of quantum
principal bundles is mainly determined by the
coalgebra structure of the quantum group. The algebra
structure enters in a few places and it is really
needed when the covariance properties are discussed
(for example we need an antipode to analyse
the transformation properties of a connection).

Let $C$ be a coalgebra and let $P$ be an algebra
and a right $C$-comodule. Assume that there
is an action $\rho : P\otimes C\otimes P \to P\otimes
C$ of $P$ on $P\otimes C$ and an element $1\in C$
such that
\roster
\item For any $u, v\in P$, $\rho(u\otimes 1 , v)
= uv\sw 0\otimes v\sw 1 $;
\item The following diagram
$$
\CD
P\otimes P @>m>>  P \\
@V{(\Delta\sb R\otimes id)}VV  @VV{\Delta\sb R}V \\
P\otimes C\otimes P @>\rho>> P\otimes C
\endCD
$$
where $m$ is a product in $P$, is commutative.
\endroster
Define $B = P\sp{coC} = \{u\in P; \Delta\sb R u = u\otimes 1\}$.
\proclaim{Lemma 3.1} $B$ is a subalgebra of $P$.
\endproclaim
\demo{Proof} Take any $u, v\in B$. Then
$$
\Delta\sb R (uv) = \rho (u\sw 0\otimes u\sw 1, v)
= \rho(u\otimes 1, v) = uv\sw 0\otimes v\sw 1 = uv\otimes 1
\qed
$$
\enddemo
\definition{Definition 3.2} Let $P$, $C$, $\rho$ and $B$
be as before.
We say that $P(B,C,\rho)$ is a {\it $C$-Galois extension}
or a {\it quantum $\rho$-principal bundle (with universal
differential structure)} if the canonical
map $\chi : P\otimes\sb B P \to P\otimes C$, $\chi: u\otimes\sb Bv
\mapsto uv\sw 0\otimes v\sw 1$ is a bijection.
\enddefinition
\example{Example 3.3} A quantum principal bundle $P(B,H)$ as
defined in \cite{BM} is a $\rho$-principal bundle with
the action $\rho: P\otimes H\otimes P\to P\otimes H$
given by $\rho(u\otimes a, v) = uv\sw 0\otimes av\sw 1$.
\endexample
\example{Example 3.4} Let $H$ be a Hopf algebra, $C$ a coalgebra and
$\pi: H\to C$ a coalgebra projection. Then $H$ is a right $C$-comodule
with a coaction $\Delta\sb R =  (id\otimes\pi)\circ\Delta$.  Denote
$1 =\pi(1)\in C$ and define $B=H\sp{coC}$ as before. Assume
that $\ker\pi$ is a minimal right ideal in $H$ such that
$\{ u - \epsilon(u) ; u\in B\} \subset \ker\pi$ (compare
Section~2). Then we can define
a canonical right action $\rho\sb 0:C\otimes H\to C$ as in the
proof of Proposition~2.5. Furthermore we define
$$
\rho(u\otimes a, v) = uv\sw 1\otimes \rho\sb 0 (a, v\sw 2),
$$
for any $u,v\in H$, $a\in C$. With these definitions $H(B,C,\rho)$ is
a quantum $\rho$-principal bundle.\endexample

\demo{Proof} First we need to show that $\rho : H\otimes C\otimes H \to
H\otimes C$ is a right action and it has the properties (1) and (2). Since
$\rho\sb 0 $ is a right action, for any $u,v,w\in H$, $a\in C$ we compute
$$\split
\rho(u\otimes a, vw)
     &=uv\sw 1w\sw 1\otimes\rho\sb 0 (a, v\sw 2 w\sw 2)
         =uv\sw 1w\sw 1\otimes \rho\sb 0(\rho\sb 0(a,v\sw 2), w\sw 2) \\
     &=\rho(uv\sw 1\otimes \rho\sb 0(a,v\sw 2), w)
         =\rho(\rho(u\otimes a,v),w),
\endsplit$$
and thus $\rho$ is an action as required. Furthermore
$$
\rho(u\otimes 1, v) = uv\sw 1\otimes \rho\sb 0(1,v\sw 2) = uv\sw 1\otimes
\pi(v\sw 2) = uv\sw 0\otimes v\sw 1,
$$
and
$$
\rho(u\sw 0\otimes u\sw 1, v) = u\sw 1v\sw 1\otimes \rho\sb 0(\pi(u\sw 2), v\sw
2)
   = u\sw 1v\sw 1\otimes \pi(u\sw 2v\sw 2) = \Delta\sb R(uv).
$$
Therefore $\rho$ has all the required properties.

To prove that the canonical map $\chi$ is bijective we first note that,
by assumption,
 $\ker\pi\subset m\circ (\ker\pi\mid\sb B\otimes H)$ and then use
a suitably modified argument of the proof of Lemma~5.2. of  \cite{BM} to
deduce that $\chi$ is a bijection. It is clear that $\chi$ is a
surjection since for any  $\sum\sb k u\sb k\otimes a\sb k \in H\otimes C$
we can choose $\sum\sb k u\sb k S{v\sb k}\sw 1 \otimes\sb B {v\sb k}\sw 2
\in H\otimes\sb B H$, where $\forall k ,$ $v\sb k \in \pi\sp{-1}(a\sb k)$,
and compute
$$\split
\chi(\sum\sb k u\sb kS{v\sb k}\sw 1 \otimes\sb B {v\sb k}\sw 2) &=
 \sum\sb k u\sb k(S{v\sb k}\sw 1) {v\sb k}\sw 2 \otimes \pi({v\sb k}\sw 3)\\
&= \sum\sb k u\sb k\otimes \pi(v\sb k) = \sum\sb k u\sb k\otimes a\sb k .
\endsplit$$
Next we compute $\ker\chi\subset H\otimes\sb B H$. Take any
$\sum\sb k u\sb k \otimes\sb B v\sb k \in \ker\chi$. Then
$\sum\sb k u\sb k {v\sb k}\sw 1 \otimes \pi({v\sb k}\sw 2) =0$.
Applying $id\otimes\epsilon$ to the last equality we then
find that $\sum\sb k u\sb k v\sb k = 0$, i.e.,
$\sum\sb k u\sb k\otimes v\sb k \in \ker m$. Any $\sum\sb i w'\sb i
\otimes w''\sb i\in\ker m$ can be written as
$\sum\sb k u\sb k S{v\sb k}\sw 1 \otimes v{\sb k}\sw 2\in H\otimes H$,
where $\forall k$, $v\sb k\in\ker\epsilon$ and $u\sb k$ are linearly
independent. Thus
$$
\chi (\sum\sb i w'\sb i\otimes\sb B w''\sb i) =
\chi(\sum\sb k u\sb k S{v\sb k}\sw 1 \otimes\sb B {v\sb k}\sw 2) =
\sum\sb k u\sb k\otimes\pi(v\sb k).
$$
If $\sum\sb i w'\sb i\otimes\sb B w''\sb i \in \ker\chi$ then
$\sum\sb k u\sb k\otimes\pi(v\sb k) = 0$, thus for all $k$,
$\pi(v\sb k) = 0$. By assumption $v\sb k = \sum\sb j b\sb k\sp j
v\sb k\sp j$, where $b\sb k \sp j \in\ker\epsilon\mid\sb B =
\ker\pi\mid\sb B$. Therefore
$$\split
\sum\sb i w'\sb i\otimes\sb B w''\sb i
        &=\sum\sb k u\sb k S{v\sb k}\sw 1 \otimes\sb B {v\sb k}\sw 2
           = \sum\sb{j,k}  u\sb k (Sv\sb k\sp j\sw 1) Sb\sb k\sp j\sw 1
           \otimes\sb B b\sb k\sp j\sw 2 v\sb k\sp j\sw 2\\
        &=\sum\sb{j, k}\epsilon(b\sb k\sp j) u\sb k Sv\sb k\sp j\sw 1
           \otimes\sb B v\sb k\sp j\sw 2 = 0
\endsplit$$
So $ker\chi = 0$, and  $\chi$ is a bijection as required.
\qed\enddemo

Therefore we have shown that an embeddable $H$ space
which is a quotient space $B=H\sp{coC}$ as described in Section~2 may
be indentified with a base manifold of the generalised
quantum principal bundle, or equivalently that $H$ is
a $C$-Galois extension of $B$.
\head 4. Manin's plane as a quantum quotient space \endhead
In this section we show that Manin's plane is a quotient
space of  the quantum general linear group $\glq$.
Recall that  Manin's plane $\cq$ is defined for any non-zero $q\in \C$ as an
associative
polynomial algebra  over $\C$ generated by $1,x, y$  subject to the
relations $xy=qyx$.  It is a quantum homogeneous space of
the quantum linear group $\glq$. $\glq$ is defined as
follows. First we consider an algebra generated
by the matrix ${\bold t} = \pmatrix \alpha &\beta \\ \gamma &\delta\endpmatrix$
and the
relations
$$
\alpha\beta =q\beta\alpha, \quad \alpha\gamma =q\gamma\alpha, \quad
\alpha\delta = \delta\alpha + (q-q\sp{-1})\beta\gamma ,\tag{4-1a}
$$
$$\beta\gamma = \gamma\beta, \quad \beta\delta = \delta\beta ,
\quad \gamma\delta=q\delta\gamma .\tag{4-1b}
$$
The quantum determinant $c=\alpha\delta - q\beta\gamma$ is central
in the algebra \thetag{4-1} thus we enlarge it with $c\sp{-1}$ and
call the resulting algebra $\glq$. The quantum linear group
$\glq$ is a Hopf algebra of a matrix group type, i.e.
$$
\Delta {\bold t} = {\bold t}\overset\cdot\to\otimes{\bold t}, \quad
\epsilon{\bold t} = 1,
\quad S{\bold t} = c\sp{-1}\pmatrix \delta &-q\sp{-1}\beta \\ -q\gamma &\alpha
\endpmatrix .
$$
The left coaction of $\glq$ on $\cq$  is given by
$$
\Delta\sb L \pmatrix x\\ y\endpmatrix  = \pmatrix \alpha &\beta \\ \gamma
&\delta
\endpmatrix
\overset\cdot\to\otimes \pmatrix x\\  y\endpmatrix .
$$
$\cq$ in not only a homogeneous space of $\glq$ but also it is
an embeddable $\glq$-space.  The linear map $\kappa :\cq\to\C$ ,
$\kappa(x\sp n y\sp m) = \delta\sb{m0}$, $m,n\in \bold Z\sb{\geq 0}$
 is a character of $\cq$.
By Proposition~2.2. it induces an algebra map $i\sb\kappa :\cq\rightarrow\glq$,
which is explicitly given by $i\sb\kappa (x) = \alpha$, $i\sb\kappa(y) =
\gamma$.
The map $i\sb\kappa$ is clearly an inclusion.
Thus the right ideal $J\sb\kappa$ is generated by $\alpha -1$ and $\gamma$.
The coalgebra $C = \glq/J\sb\kappa$ may be easily computed. It is
spanned by $a\sb{m,n}=\pi(\beta\sp mc\sp n)$, $m\in {\bold Z}\sb{>0}$, $n\in
\bold Z$
and $a\sb{0,0} = 1 = \pi(1)$, where $\pi :\glq\to C$ is a canonical surjection.
To see that the $a\sb{m,n}$ really span $C$ we note that since $J\sb\kappa$
is generated by $\alpha -1$ and $\gamma$ as
a right ideal in $\glq$, every $\alpha$ which multiplies
any element of $\glq$ from the left
 is replaced by $1$ and similarly any $\gamma$
is replaced by $0$ when the resulting
element of $\glq$ is acted upon by $\pi$. Then we compute
$$\split
\pi(\alpha\sp k\beta\sp l\gamma\sp m\delta\sp n c\sp r)
  &= \pi(\beta\sp l\gamma\sp m\delta\sp n c\sp r)  \\
  &= \delta\sb{m0}\pi(\beta\sp l\delta\sp n c\sp r)  = \delta\sb{m0}q\sp{ln}
      \pi(\delta\sp n\beta\sp l c\sp r) \\
   &= \delta\sb{m0}q\sp{ln}(\pi(\alpha\delta\delta\sp{n-1}\beta\sp l c\sp r)
      - q\sp{-1}\pi(\gamma\beta\delta\sp{n-1}\beta\sp lc\sp r)) \\
   &=\delta\sb{m0}q\sp{ln}\pi(\delta\sp{n-1}\beta\sp lc\sp{r+1}) =\ldots
      =\delta\sb{m0}q\sp{ln}a\sb{l,r+n}.
\endsplit$$
Therefore any element of $C=\pi(\glq)$ may be expressed as a linear
combination of $a\sb{m,n}$.

The coalgebra structure of $C$ is found from the coalgebra structure
of $\glq$, since $\Delta\sb C = (\pi\otimes\pi)\circ\Delta\sb{\glq}$.
Explicitly
$$
\Delta a\sb{m,n} = \sum\sb{k=0}\sp m {\binom mk}\sb q a\sb{k,n}
\otimes a\sb{m-k, n+k}, \quad \epsilon(a\sb{m,n}) = \delta\sb{m0},
\tag{4-2}
$$
where the quantum binomial coefficients are defined by
$$
{\binom mk}\sb q  = \frac{[m]\sb q!}{[m-k]\sb q! [k]\sb q!}, \quad
[m]\sb q = \frac{q\sp m - q\sp{-m}}{q-q\sp{-1}}, \quad
[m]\sb q ! = \prod\sb{k=1}\sp m [k]\sb q, \quad [0]\sb q! = 1.
$$

The next step in the identification of $\cq$ as a quantum quotient
space consists of computing the fixed point subalgebra
$B = \glq\sp{coC}$. For a general monomial
$\alpha\sp k\gamma\sp l\beta\sp m\delta\sp n c\sp r \in \glq$,
$k, l, m, n \in \bold Z\sb{\ge 0}$, $r\in\bold Z$, we find
$$\split
\Delta\sb R(\alpha\sp k\gamma\sp l\beta\sp m\delta\sp n c\sp r)
   &= \alpha\sp k\gamma\sp l c\sp r \sum\sb{i=0}\sp m\sum\sb{j=0}
     \sp n q\sp{j(m-i)} {\binom mi}\sb q {\binom nj}\sb q \\
   & \times\alpha\sp{m-i}\beta\sp i\gamma\sp{n-j}\delta\sp j
     \otimes a\sb{m+n -(i+j), i+j+r}.\endsplit\tag{4-3}
$$
The right hand side of \thetag{4-3} has the form $u\otimes 1$
for some $u\in\glq$ if and  only if $m=n=r=0$. Thus $B$
is spanned by all $\alpha\sp k\gamma\sp l$. Therefore
$B\subset i\sb\kappa(\cq)$ and since $i\sb\kappa(\cq)\subset B$
by Proposition~2.5 we conclude that $\cq \cong \glq\sp{coC}$.
By Example~3.4 $\glq(\cq, C, \rho)$ is a quantum principal
$\rho$-bundle. The action $\rho\sb 0:C\otimes\glq\to C$ is
given explicitly by
$$
\rho\sb 0(a\sb{i,j},\alpha\sp k\beta\sp l\gamma\sp m\delta\sp n c\sp r)
= \delta\sb{m0}q\sp{i(n-k)+ln}a\sb{i+l,j+n+r}.
$$

We can now proceed to define an algebra structure
on $C$ so that it becomes a Hopf algebra. We define
the product in $C$ by
$$
a\sb{k,l}a\sb{m,n} = q\sp{lm-kn}a\sb{k+m , l+n}.
$$
First we notice that $a\sb{0,0} =1$ is the unit element
with respect to this product. Next we show that this
product is compatible with the coalgebra structure
of $C$. We compute
$$\split
\Delta(a\sb{k,l})\Delta(a\sb{m,n})
  &=\sum\sb{i=0}\sp k\sum\sb{j=0}
     \sp m {\binom ki}\sb q {\binom mj}\sb q
     a\sb{i,l}a\sb{j,n}\otimes a\sb{k-i, l+i}
     a\sb{m-j, n+j} \\
  &=q\sp{lm-nk} \sum\sb{i=0}\sp k\sum\sb{j=0}
     \sp n q\sp{im-kj} {\binom ki}\sb q {\binom mj}\sb q
     a\sb{i+j , l+n}\otimes a\sb{k+m -(i+j), l+n+i+j}\\
  &=q\sp{lm-nk} \sum\sb{r=0}\sp{k+m}{\binom {k+m}i}\sb q
     a\sb{r, l+n}\otimes a\sb{k+m-r , l+n+r} \\
  &= q\sp{lm-kn}\Delta(a\sb{k+m , l+n} )
     = \Delta(a\sb{k,l}a\sb{m,n}).
\endsplit$$
The third equality is a consequence of the
following property of the $q$-deformed binomial
coefficients
$$
\forall r\in [0, k+m], \quad \sum\sb{i=0}\sp k\sum\sb{j=0}
     \sp n q\sp{im-kj} {\binom ki}\sb q {\binom mj}\sb q
= {\binom {k+m}r}\sb q.
$$
Clearly the counit of $C$ is an algebra homomorphism.
Before we define an antipode we show that $C$ is
a polynomial algebra. Let $a=a\sb{0,1}$, $a\sp{-1} = a\sb{0,-1}$,
 $b = a\sb{1,0}$.
Then for any $m\in\bold Z\sb{\geq 0}$, $n\in \bold Z$,
$$
a\sb{m,n} = q\sp{-mn}a\sp n b\sp m, \qquad ab=q\sp{2}ba, \qquad aa\sp{-1}
= a\sp{-1} a = 1.
$$
Therefore $C$ is a polynomial algebra indeed, and it is isomorphic
to $\bold C\sb{q\sp 2}\sp 2[x\sp{-1}]$. The coalgebra structure
of $C$ written in terms of $a$ and $b$ reads
$$
\Delta a\sp{\pm 1}= a\sp{\pm 1}\otimes a\sp{\pm 1}, \qquad \Delta b = 1\otimes
b
+ b\otimes a, \qquad \epsilon(a\sp{\pm 1}) =1, \qquad \epsilon(b) =0
$$
and hence the antipode is defined as $S a\sp{\pm 1}=a\sp{\mp 1}$,
$S b = -ba\sp{-1}$.

We have just shown that $C$ may be equipped with
an algebra structure of $\bold C\sb{q\sp 2}\sp 2[x\sp{-1}]$,
and then the coalgebra structure of $C$ becomes a
standard coalgebra structure of the latter. Therefore
we have proven
\proclaim{Theorem 4.1}
$$
\cq = \glq\sp{co\bold C\sb{q\sp 2}\sp 2[x\sp{-1}]}.
$$
\endproclaim

Notice that clearly neither $\pi: \glq\to\bold C\sb{q\sp 2}
\sp 2[x\sp{-1}]$ nor $\Delta\sb R = (id\otimes \pi)\circ\Delta :
\glq\to\glq\otimes\bold C\sb{q\sp 2}\sp 2[x\sp{-1}]$ are algebra
maps. Still, following the proposal of Section~3
we can analyse the generalised principal bundle
$\glq(\cq, \bold C\sb{q\sp 2}\sp 2[x\sp{-1}], \rho, \pi)$. In
particular we can truly develop a gauge theory,
define connections and their curvature, following closely
the quantum group gauge theory introduced in \cite{BM}.

\head 5. Podle\'s' sphere as a quantum quotient space \endhead
In this section we prove that the quantum two-sphere
is a quantum quotient space in the sense expalined in
Section~2. In our presentation of
the quantum sphere we follow the conventions of \cite{NM}.

The general quantum two-sphere $S\sb q \sp 2(\mu ,\nu)$ is
a polynomial algebra generated by the unit and $x$, $y$, $z$,
and the relations
$$\alignat 2
 xz & = q\sp 2zx, & \qquad xy & = -q(\mu -z)(\nu +z), \\
 yz & = q\sp{-2}zy, & \qquad yx & = -q(\mu -q\sp{-2}z)(\nu +q\sp{-2} z),
\endalignat$$
where $\mu, \nu $ and $q\neq 0$ are real parameters, $\mu\nu \geq 0$,
$(\mu ,\nu) \neq (0,0)$.
The quantum sphere is a $*$-algebra with the
$*$-structure $x\sp * = -q y$, $z\sp * =z$.

The quantum sphere $S\sp 2\sb q(\mu ,\nu)$ is an $SU\sb q(2)$
homogeneous quantum space. $SU\sb q(2)$ is defined as a quotient
of $\glq$ by the relation $c =1$, and has a $*$-structure given
by $\delta = \alpha\sp *$, $\gamma = -q\sp{-1} \beta\sp *$. The
coaction of $SU\sb q(2)$ on $\sq$ is defined as follows.
Let $\phi\sb - =x$, $\phi\sb 0 =(1+q\sp{-2})\sp{-1/2}(\mu -\nu
- (1+q\sp{-2}) z)$, $\phi\sb + =y$. Then
$$
\Delta\sb L \pmatrix \phi\sb -\\ \phi\sb 0 \\ \phi\sb +\endpmatrix
= \pmatrix \alpha\sp 2 & (1+q\sp{-2})\sp{1/2}\alpha\beta &
\beta\sp 2 \\ (1+q\sp{-2})\sp{1/2}\alpha\gamma & 1 +
(q +q\sp{-1})\beta\gamma & (1+q\sp{-2})\sp{1/2}\beta\delta \\
\gamma\sp 2 & (1+q\sp{-2})\sp{1/2}\gamma\delta & \delta\sp 2
\endpmatrix \overset\cdot\to\otimes \pmatrix \phi\sb -\\
\phi\sb 0 \\ \phi\sb +
\endpmatrix .
$$
The quantum sphere $\sq$ is not only a quantum homogeneous
space but also it is an embeddable $SU\sb q(2)$-space.
There is a $*$-character $\kappa : \sq\to\C$ given by
$$
\kappa(x) = q\sqrt{\mu\nu}, \qquad \kappa(y) = -\sqrt{\mu\nu},
\qquad \kappa(z) = 0 .
$$
Therefore there is also a $*$-algebra homomorphism
$i\sb\kappa : \sq\rightarrow SU\sb q(2)$, which reads
explicitly
$$
i\sb\kappa(x) = \sqrt{\mu\nu}(q\alpha\sp 2 - \beta\sp 2) +
(\mu -\nu)\alpha\beta,
$$
$$
i\sb\kappa(y) = \sqrt{\mu\nu}(q\gamma\sp 2 - \delta\sp 2) +
(\mu -\nu)\gamma\delta,
$$
$$
i\sb\kappa(z) = -\sqrt{\mu\nu}(q\alpha\gamma - \beta\delta) -
(\mu -\nu)\beta\gamma ,
$$
and is clearly an inclusion. From now on we assume that
$\mu\neq\nu$ (but see also Remark~5.5). In this case
$\sq$ depends on two real parameters only, namely $q$
and $p = \frac{\sqrt{\mu\nu}}{\mu-\nu}$.
By Proposition~2.4. the inclusion $i\sb\kappa$ induces a
coideal $J\sb\kappa\subset SU\sb q(2)$, generated as a right ideal in
$SU\sb q(2)$ by the following three elements
$$
p(q\alpha\sp 2 - \beta\sp 2) + \alpha\beta -pq, \quad
p(q\gamma\sp 2 - \delta\sp 2) + \gamma\delta + p, \quad
p(q\alpha\gamma - \beta\delta) + q\beta\gamma .
$$
Therefore we can construct the coalgebra
$C(p) = SU\sb q(2)/J\sb \kappa$, and the corresponding
quotient space $B(p) = SU\sb q(2)\sp{coC(p)}$  as
described in Section~2. At the end of this procedure
we identify $B(p)$ with $\sq$, $\mu\neq \nu$.
We start with the coalgebra
$C(p)$.

\proclaim{Proposition 5.1}
$C(p)$ is a vector space spanned by $1=\pi(1)$, $x\sb n=
\pi(\alpha\sp n)$ and $y\sb n = \pi(\delta\sp n)$, where
$\pi :SU\sb q(2)\to C(p)$ is a canonical surjection and $n\in\bold Z\sb{>0}$.
\endproclaim
\demo{Proof} For any $u\in SU\sb q(2)$ we use the explicit
form of the generators
of $J\sb\kappa$, and the relations in $SU\sb q(2)$ to find that
$$\split
\pi(\beta u)
  &=\pi(\beta\alpha\delta u) - q\pi(\beta\gamma\beta u)\\
  &=q\sp{-1} \pi (\alpha\beta\delta u) - q\pi(\beta\gamma\beta u)\\
  &=-p\pi(\alpha\sp{2}\delta u) +pq\sp{-1}\pi (\beta\sp 2\delta u)\\
  &+p\pi(\delta u) +pq\pi(\alpha\gamma\beta u) - p\pi(\beta\delta
   \beta u) \\
  &=p\pi(\delta u) - p\pi(\alpha u),
\endsplit\tag{5-1a}$$
and similarly
$$
\pi(\gamma u) = p\pi(\delta u) - p\pi(\alpha u). \tag{5-1b}
$$
{}From \thetag{5-1} it follows that
for any $u\in SU\sb q(2)$, $\pi(u\beta\sp m\gamma\sp n)
= \pi(u\beta\sp{m+n})$.
Since $SU\sb q(2)$ is spanned by the monomials
$\alpha\sp m\beta\sp k\gamma\sp l$, $\delta\sp m\beta\sp k
\gamma\sp l$ (cf. Lemma~7.1.2 of \cite{CP}) it suffices to prove that the
following elements of $C(p)$,
$$
 a\sb{k-}\sp{(n)} = \pi(\delta\sp k\beta\sp{n-k}),\qquad
a\sb{k+}\sp{(n)} = \pi(\alpha\sp k\beta\sp{n-k}),
\tag{5-2}$$
where $n\in\bold Z\sb{>0}$, $k=0,1,\ldots , n$,
can be expressed as linear combinations of $1$, $x\sb m$,
$y\sb m$. Clearly $a\sp{(n)}\sb{0 -} = a\sp{(n)}\sb{0 +}$.
Thus we simply write $a\sp{(n)}\sb{0}$. Also, $a\sp{(n)}\sb{n+}
 = x\sb n$ and $a\sp{(n)}\sb{n-}
 = y\sb n$. For $n=1$,
$a\sp{(1)}\sb{0} = \pi(\beta) = p(y\sb 1 - x\sb 1)$.
For a general $n$ we apply the rules \thetag{5-1}
to $a\sp{(n)}\sb{k\pm}$ and we express the latter
in terms of $a\sp{(m)}\sb{l\pm}$, $m<n$ and $x\sb n$, $y\sb n$.
We make the inductive assumption that for all $m<n$,
$a\sp{(m)}\sb{l\pm}$ can be written as linear combinations of $1$,
$x\sb r$, $y\sb r$. Therefore, for $n\geq 2$ we arrive at the
system of equations
$$\aligned
a\sp{(n)}\sb{k\pm} &\pm pq\sp{\pm k}a\sp{(n)}\sb{k+1\pm}
     \mp pq\sp{\mp(k-1)}a\sp{(n)}\sb{k-1\pm} = \pm
     pq\sp{\pm k}a\sp{(n-2)}\sb{k-1\pm}, \\
a\sp{(n)}\sb{0} &- pa\sp{(n)}\sb{1 -}
     +pa\sp{(n)}\sb{1 +} = 0,
\endaligned\tag{5-3}$$
where $k = 1,2, \ldots , n-1$. This is a system of
$2n-1$ equations with $2n-1$ unknowns provided
that the right hand sides and $x\sb n$, $y\sb n$
are treated as known parameters. Obviously it has a solution
if its determinant is non-zero. The determinant $D\sb n$ of
the system \thetag{5-3} may be easily computed.
It does not depend on $q$ and, by the Laplace theorem,
 it can be reduced to
the determinant of the following $2n-1\times 2n-1$ matrix
$$
\pmatrix
1&-p&p&0&0&0&\hdots&0&0&0&0&0\\
p&1&0&-p&0&0&\hdots&0&0&0&0&0\\
-p&0&1&0&p&0&\hdots&0&0&0&0&0\\
0&p&0&1&0&-p&\hdots&0&0&0&0&0\\
\hdotsfor{12}\\
0&0&0&0&0&0&\hdots&0&1&0&-p&0\\
0&0&0&0&0&0&\hdots&-p&0&1&0&p\\
0&0&0&0&0&0&\hdots&0&p&0&1&0\\
0&0&0&0&0&0&\hdots&0&0&-p&0&1
\endpmatrix
\tag{5-4}
$$
Again by the Laplace theorem $D\sb n$ can be further developed
to give
$$
D\sb n = A\sb{2n-2} + p\sp 2(A\sb{2n-3} +A\sb{2n-4}) +
p\sp 4 A\sb{2n-5},
$$
where $A\sb m$ is zero for negative $m$, $A\sb 0 =1$
and for any $m < 2n-1$, $A\sb m$ is the determinant
of the matrix obtained from \thetag{5-4} by removing
first $2n-1-m$ rows and columns. The determinants $A\sb m$
are the standard ones and we finally obtain the determinant
of the system \thetag{5-3} as a polynomial
$$
D\sb n = P\sb{n-1}(p\sp 2)
  \equiv \sum\sb{k=0}\sp{n-1} \binom{2n-1-k}{k} p\sp{2k}.
$$
For any $x\in \bold R\sb{\geq 0}$, $P\sb n(x)\geq 1$,
and hence $D\sb n \neq 0$ for any real $p$. Therefore the system
\thetag{5-3} always has a solution and the coalgebra $C(p)$
is spanned by $x\sb n$, $y\sb n$, $n\in\bold Z\sb{>0}$
and $1$ as required. \qed\enddemo

The vector space $C(p)$ has a coalgebra structure
induced by $\pi$ from the coalgebra structure of
$SU\sb q(2)$. The coproduct reads explicitly
$$
\Delta x\sb n = \sum\sb{k=0}\sp{n-1} q\sp{-(n-k)k}
{\binom nk}\sb q a\sp{(n)}\sb{k+}\otimes a\sp{(n)}\sb{k+}, \qquad
\Delta y\sb n = \sum\sb{k=0}\sp{n-1} q\sp{(n-k)k}
{\binom nk}\sb q a\sp{(n)}\sb{k-}\otimes a\sp{(n)}\sb{k-} ,
$$
where $a\sb{k\pm}\sp{(n)}$ are given by \thetag{5-2}.
Therefore the coalgebra $C(p)$ is cocommutative.
\remark{Remark 5.2} It is
an interesting problem, whether it is possible to
define a Hopf algebra structure on $C(p)$.
For example,  for $n = 1$ we have
$$\align
\Delta x\sb 1 &= (1+p\sp 2)x\sb 1\otimes x\sb 1
- p\sp 2(x\sb 1\otimes y\sb 1 +y\sb 1\otimes x\sb 1
- y\sb 1\otimes y\sb 1), \\
\Delta y\sb 1 &= (1+p\sp 2)y\sb 1\otimes y\sb 1
- p\sp 2(x\sb 1\otimes y\sb 1 +y\sb 1\otimes x\sb 1
- x\sb 1\otimes x\sb 1).
\endalign$$
If we define
$$
x'\sb 1 = \frac{1}{\mu -\nu}(\mu x\sb 1 -\nu y\sb 1),
\qquad y'\sb 1 = \frac{1}{\mu -\nu}(\mu y\sb 1 -\nu x\sb 1)
$$
then $x'\sb 1$ and $y'\sb 1$ are group-like, i.e.
$\Delta x'\sb 1 = x'\sb 1\otimes x'\sb 1$ and
$\Delta y'\sb 1 = y'\sb 1\otimes y'\sb 1$.
If it were posible to
define a new basis of $C(p)$ consisting only of group-like elements
then clearly we would be able to solve the above problem and
make $C(p)$ into a Hopf algebra of functions
on $U(1)$.
\endremark
\remark{Remark 5.3}
According to \cite{P1} quantum spheres can also
be defined for a discrete series of complex numbers $p$
given by $p\sp 2 = - (q\sp k +q\sp{-k})\sp{-2}$, $k=1,2,\ldots$.
It is shown in \cite{P2} that such quantum spheres are
$*$-embeddable in $SU\sb q(2)$ for $k=1$.

One easily finds that $P\sb n(x-1/4) =
\sum\sb{k=0}\sp{n} c\sp n\sb kx\sp k$, where
$$
c\sp n\sb k = \sum\sb{l=k}\sp n(-1/4)\sp{l-k}\binom{2n+1-l}{l}
\binom lk
$$
 For any $n$ and any $0\leq k\leq n$, $c\sp n\sb k
\geq c\sp n\sb 0 = (n+1)/4\sp n$ and thus all the
coefficients $c\sp n\sb k $ are positive. Therefore
$P\sb n(x-1/4)\neq 0$ for all real $x\geq 0$. This
implies that the determinants $D\sb n$ of the proof
of Proposition~5.1 are non-zero provided that
$p\sp 2\geq -1/4$. Since for any $q$, $q +q\sp{-1}\geq 2$ we
see that the
assertion of Proposition~5.1 holds for the exceptional
 quantum spheres too.
\endremark
\proclaim{Propositon 5.4}
Let $C(p)$
 be a coalgebra described in Proposition~5.1
and let $B(p) = SU\sb q(2)\sp{coC(p)}$.
Then $i\sb\kappa(\sq)
= B(p)$ for all $\mu\neq \nu$  such that
 $p= \frac{\sqrt{\mu\nu}}{\mu - \nu}$.
\endproclaim
\demo{Proof} By Proposition~2.5, $i\sb\kappa(\sq)\subset B(p)$,
therefore
we need to show that $B(p) \subset i\sb\kappa(\sq)$. We introduce
the grading $d: SU\sb q(2)\to \bold Z$ by
$$
d(\alpha) = d(\beta) = 1, \quad d(1) = 0, \quad d(\gamma) =
d(\delta) = -1, \quad d(uv) = d(u) + d(v),
\tag{5-5}$$
for any monomials $u, v\in SU\sb q(2)$. A set of all elements of
$SU\sb q(2)$ of
degree $k\in\bold Z$ forms a vector subspace of $SU\sb q(2)$,
which we denote by $SU\sb q(2)\sp{(k)}$, and $SU\sb q(2)
= \bigoplus\sb{k\in\bold Z} SU\sb q(2)\sp{(k)}$. Moreover if
$\Delta u = \sum\sb i u'\sb i \otimes u''\sb i $ for any
$u\in SU\sb q(2)\sp{(k)}$, then for all $i$, $d(u'\sb i) = k$.
To see that the last statement is true we can explicitly verify it
for $\alpha$, $\beta$, $\gamma$, $\delta$ and then use definition
\thetag{5-5} of $d$ to prove it for any $SU\sb q(2)$. Therefore
$d$ induces a grading of $B(p)$ and $B(p) = \bigoplus\sb{k\in\bold
Z}B(p)\sp{(k)}$.

Next we notice that $B(p)$ is contained in the subalgebra of $SU\sb q(2)$
spanned by monomials of even degree. Therefore for any $k\in\bold Z$,
$B(p)\sp{(2k+1)} =0$.

To prove the required inclusion we observe that due to
the form of $\pi$ and $C(p)$, $B(p)$ is a deformation
of $B(0)$, i.e., $B(0) = \underset{p\to 0}\to\lim B(p)$. We denote
by $B(p)\sb{2n}\sp{(2k)}$ the vector space of homogeneous polynomials
$u\in B(p)$ of degree $2n$ such that $d(u)=2k$, $|k|\leq n$.
Notice that $B(p)\sp{(2k)}\sb{2n}$ and $B(p)\sp{(2k)}\sb{2l}$ need
not be distinct for $l\neq n$.
$B(0)\sb{2n}\sp{(2k)}$
is spanned by $\alpha\sp{m}\beta\sp{n+k -m}\gamma\sp{n-m}
\delta\sp{m-k}$,
where $m = k, k+1, \ldots , n$ for $k\geq 0$ and $m = 0, 1,\ldots ,
n+k$ for $k<0$, and hence is $n - |k|+1$-dimensional. This is
exactly the dimension of $i\sb\kappa(\sq)\sb{2n}\sp{(2k)}$.
Suppose that $B(p)\sb{2n}\sp{(2k)}$ is at least $n - |k|+2$
dimensional. Then we can find $u\in B(p)\sb{2n}\sp{(2k)}$ that
does not contain any of the monomials spanning $B(0)\sb{2n}\sp{(2k)}$.
If $\underset{p\to 0}\to\lim u \neq 0$, then we would obtain
that $B(0)\sb{2n}\sp{(2k)}$ is at least $n-|k|+2$ dimensional,
hence contradiction. By $\underset\sb{p\to 0}\to\lim u$ here we mean
the polynomial obtained from $u$ by replacing its coefficients
with their $p=0$ limits. Assume that $\underset\sb{p\to 0}\to\lim u =0$.
The polynomial $u$ may be written as a linear combination of
monomials of degree $2n$ with coefficients that vanish
as polynomials when $p$ tends to $0$. Therefore there
exists a positive integer $m$ such that $\underset\sb{p\to 0}\to
\lim p\sp{-m}u$
exists, is finite and non-zero, and is an element of
$B(0)\sb{2n}\sp{(2k)}$. Thus we have a contradiction again.
Since the above argument does not depend on $n$ and $k$,
and $i\sb\kappa(\sq)\subset B(p)$ we conclude that
$i\sb\kappa(\sq) = B(p)$. \qed \enddemo

Therefore we have shown that for $\mu\neq \nu$ the quantum
sphere $\sq$ is a quantum quotient space. By Example~3.4
we also have a principal $\rho$-bundle, $SU\sb q(2)(\sq, C(p),
\rho, \pi)$.

\remark{Remark 5.5} When $\mu=\nu\neq 0$ the coideal $J\sb\kappa$
is generated as a right ideal in $SU\sb q(2)$
by the following elements
$$
q\alpha\sp{2} - \beta\sp 2 -q,\qquad q\gamma\sp 2 -\delta\sp 2 +1,
\qquad q\alpha\gamma - \beta\delta .
$$
Therefore for any $u\in SU\sb q(2)$,
$$
\pi(\delta u) = \pi(\alpha u), \qquad \pi(\beta u ) =\pi(\gamma u),
\qquad \pi(\gamma\sp 2 u ) = q\pi(\alpha\sp 2 u) - q\pi(u),
$$
and hence the coalgebra $C = SU\sb q(2)/J\sb\kappa$ is spanned by
$1=\pi(1)$, $x\sb n = \pi(\alpha\sp n)$, $y\sb n =
\pi(\alpha\sp{n-1}\gamma)$, $n\in\bold Z\sb{\geq 1}$. We conjecture
that also for this case
$S\sp 2\sb q(\mu,\mu)\cong SU\sb q(2)\sp{coC}$.
\endremark

\head 6. Conclusions \endhead

In this paper we have shown that certain embeddable
quantum homogeneous spaces may be viewed as quantum
quotient spaces. The examples of such quantum
embeddable spaces include
the general quantum two-sphere $\sq$ and
the quantum plane $\cq$. The interpretation of
quantum embeddable spaces presented in
this paper seems specially interseting
from the point of view of quantum group gauge theory,
the suitable generalisation of which we have also
proposed. We think that it would be intersting
and indeed desirable to develop further this generalisation
of quantum group gauge theory, and in particular,
to construct connections on the
quantum spaces described in this paper. For example
this would allow
for extending the construction of the Dirac $q$-monopole
of \cite{BM} to general quantum spheres.

\subhead Acknowledgement\endsubhead
I would like to thank Shahn Majid for suggesting to me the
proofs of Propositions 2.4 and 2.5(2).


\baselineskip 16pt
\Refs
\widestnumber\key{MM}
\ref\key BM \by T. Brzezi\'nski and S. Majid
\paper Quantum Group Gauge Theory on Quantum Spaces
\jour Commun. Math. Phys. \yr1993 \vol 157 \page 591
\moreref \jour ibidem. \vol 167 \yr1995
\page 235 \lang Erratum \endref
\ref\key CP \by V. Chari and A. Pressley \book
A Guide to Quantum Groups \publ Cambridge University
Press \yr 1994 \endref
\ref\key C \by A. Connes \book A lecture given at the
Conference on Non-commutative Geometry and Its
Applications \publ Castle T\v re\v s\v t. Czech Republic,
May 1995\endref
\ref\key M \by Y.I. Manin \book Quantum Groups and
Non-commutative Geometry \publ Montreal Notes \yr 1989
\endref
\ref\key NM \by M. Nuomi and K. Mimachi \jour Commun.
Math. Phys. \vol 128 \yr 1990 \page 521 \paper Quantum
2-spheres and Big $q$-Jacobi Polynomials\endref
\ref\key P1 \by P. Podle\'s \paper Quantum Spheres
\jour Lett. Math. Phys. \vol 14 \yr 1987 \page 193
\endref
\ref\key P2 \bysame \paper Symmetries of Quantum
Spaces. Subgroups and Quotient Spaces of Quantum
$SU(2)$ and $SO(3)$ Groups \jour Commun. Math. Phys.
\vol 170 \yr 1995 \page 1 \endref
\ref \key S \by H.-J. Schneider \paper Principal
homogeneous spaces for arbitrary Hopf algebras \jour
Israel J. Math \vol 72 \yr 1990 \page 167
\moreref \paper Representation theory of Hopf-Galois
extensions \jour ibidem. \vol 72 \yr 1990 \page 196
\endref
\ref\key Sw \by M.E. Sweedler \book Hopf algebras
\publ Benjamin \yr 1969 \endref
\endRefs
\enddocument